\documentclass[final,12pt]{elsarticle}

\usepackage{amssymb}

\usepackage{braket,amsfonts}
\usepackage{mathtools}
\usepackage[amsmath]{ntheorem}
\usepackage{array}
\usepackage{multirow}
\usepackage{enumitem}
\usepackage{xspace}
\usepackage{arydshln}
\usepackage{lineno}
\usepackage{subfig}
\usepackage{algcompatible}
\usepackage[linesnumbered,lined]{algorithm2e}
\RestyleAlgo{ruled}
\usepackage{hyperref}
\usepackage[capitalise,nameinlink,sort]{cleveref}

\usepackage{xcolor}
\usepackage{bm}
\newcommand{\norm}[1]{\lVert#1\rVert}
\newcommand{\R}{\mathbb R}
\newcommand{\N}{\mathbb N}

\newcommand{\be}{\mathbf e}
\newcommand{\bff}{\mathbf f}

\newcommand{\bp}{\mathbf p}
\newcommand{\br}{\mathbf r}
\newcommand{\bs}{\mathbf s}
\newcommand{\bu}{\mathbf u}
\newcommand{\bv}{\mathbf v}
\newcommand{\bx}{\mathbf x}
\newcommand{\by}{\mathbf y}

\newcommand{\calo}{\mathcal O}

\newcommand{\wh}{\widehat}

\DeclareMathOperator*{\argmax}{\arg\!\max}

\DeclareMathOperator{\Range}{Range}

\DeclareMathSymbol{\mrq}{\mathord}{operators}{`'}
\DeclareMathSymbol{\mlq}{\mathord}{operators}{``}

\newcommand{\BlackBox}{\rule{1.5ex}{1.5ex}}  
\ifdefined\proof
    \renewenvironment{proof}{\par\noindent{\bf Proof\ }}{\hfill\BlackBox\\[2mm]}
\else
    \newenvironment{proof}{\par\noindent{\bf Proof\ }}{\hfill\BlackBox\\[2mm]}
\fi
\newtheorem{theorem}{Theorem}[section]

 \newtheorem{remark}[theorem]{Remark}

\newtheorem{experiment}[theorem]{Experiment}

\crefname{experiment}{Experiment}{Experiments}
\crefname{proposition}{Proposition}{Propositions}
\crefname{definition}{Definition}{Definitions}
\crefname{theorem}{Theorem}{Theorems}
\crefname{lemma}{Lemma}{Lemmas}
\usepackage{pgfplots}
 \pgfplotsset{compat=newest,legend style={font=\scriptsize,row sep=-0.05cm,/tikz/every odd column/.append style={column sep=0.01cm}}} 
 \usepgfplotslibrary{groupplots}
  \newcommand{\exponent}{2}
\newlength\figureheight
\newlength\figurewidth

\newcommand{\smtxa}[2]{
{\mbox{\scriptsize
$\left[\!\!
\begin{array}{#1}
#2
\end{array} \!\! \right]$}}}

\newcommand{\mtxa}[2]{
\left[
\begin{array}{#1}
#2
\end{array}
\right]}

\usepackage{color}
\definecolor{orange}{rgb}{1.0,0.4,0}
\newcommand{\pg}[1]{{\color{black}#1}}

\journal{xxxxx}

\begin{document}

\begin{frontmatter}


\author{Perfect~Y.~Gidisu\corref{cor1}}
\ead{p.gidisu@tue.nl}
\cortext[cor1]{Corresponding author}

\title{A DEIM-CUR factorization with iterative SVDs\tnoteref{label1}}

\author{Michiel E.~Hochstenbach}
\ead{m.e.hochstenbach@tue.nl}


\begin{abstract}
A CUR factorization is often utilized as a substitute for the singular value decomposition (SVD),  especially when a concrete interpretation of the singular vectors is challenging. Moreover, if the original data matrix possesses properties like nonnegativity and sparsity, a CUR decomposition can better preserve them compared to the SVD. An essential aspect of this approach is the methodology used for selecting a subset of columns and rows from the original matrix. This study investigates the effectiveness of \emph{one-round sampling} and iterative subselection techniques and introduces new iterative subselection strategies based on iterative SVDs. One provably appropriate technique for index selection in constructing a CUR factorization is the discrete empirical interpolation method (DEIM). Our contribution aims to improve the approximation quality of the DEIM scheme by iteratively invoking it in several rounds, in the sense that we select subsequent columns and rows based on the previously selected ones. Thus, we modify $A$ after each iteration by removing the information that has been captured by the previously selected columns and rows. We also discuss how iterative procedures for computing a few singular vectors of large data matrices can be integrated with the new iterative subselection strategies. We present the results of numerical experiments, providing a comparison of one-round sampling and iterative subselection techniques, and demonstrating the improved approximation quality associated with using the latter.
\end{abstract}

\end{frontmatter}

\section{Introduction}
\label{sec:intro}

In data analysis applications and machine learning, the data set is often represented by a matrix $A\in \R^{m\times n}$, which is usually large. In many cases, a key step in the analysis is to approximate the data using a few features and/or a few data points so that one can easily manipulate, understand, and interpret the data. The optimal approximation is obtained by the truncated singular value decomposition (TSVD). On the other hand, this approximation problem can also be reduced to identifying a good subset of columns and rows in the data matrix, a CUR factorization. A CUR decomposition is an alternative solution to the SVD, motivated by the fact that in several applications, it can be challenging to have a concrete interpretation of the singular vectors. Additionally, the singular vectors fail to preserve properties such as nonnegativity and sparsity, if the original data matrix has these. 
Theoretical computer science and numerical linear algebra communities have extensively studied column subset selection (index selection) algorithms. In the latter,  researchers have primarily focused on deterministic algorithms, which exploit the SVD or rank-revealing QR factorizations that select columns by pivoting rules \cite{bischof1991structure,chandrasekaran1994,gu1996,Voronin,Goreinov2010}.

 On the other hand,  researchers in the theoretical computer science community have predominantly directed their attention towards selecting ``optimal" columns using randomized algorithms with provable error bounds \cite{optimalboutsidis,deshpande2006matrix,Drineas,frieze2004fast,guruswami2012optimal}. { Randomized algorithms sometimes offer reduced computational complexity compared to deterministic methods. Nevertheless, these techniques usually necessitate oversampling columns and rows beyond the desired target rank $k$ to achieve strong provable approximation guarantees. 
 
 In \cite{hamm2020stability}, the authors explore various random sampling methods and analyze their computational complexities, along with assessing the stability of the methods when subjected to perturbations in both the probabilities and the underlying matrix. Commonly employed sampling distributions in these algorithms include the uniform distribution (${ \sf pr}_j = 1/n$, where $n$ represents the number of columns in $A$) \cite{Chiu}, squared-norm distribution (${ \sf pr}_j = \|\textbf{A}(:, j)\|^2 / \|\textbf{A}\|_F^2$) \cite{frieze2004fast}, and leverage scores distribution (${ \sf pr}_j = \frac{1}{k} \|\textbf{V}_k(j, :)\|^2$, where $\textbf{V}_k$ contains the $k$-leading right singular vectors) \cite{Mahoney,Drineas}. We denote the spectral (2-norm) and the Frobenius norm by $\norm{\cdot}$ and $\norm{\cdot}_F$, respectively. 
 
Deterministic algorithms based on a derandomization of the volume
sampling algorithm \cite{deshpande2006matrix} have been proposed in \cite{cortinovis2020low,deshpande2010efficient}. In \cite{boutsidis2008selecting,dong2023simpler,dong2024robust}, algorithms leveraging the advantages of both randomized and deterministic methods are introduced.} This paper focuses on the deterministic algorithms that exploit the SVD for index selection; in particular, the discrete empirical interpolation method (DEIM) \cite{Barrault,Chaturantabut,Drmac,Sorensen}.

The notations $C^+$ and $C^\top$ denote the Moore--Penrose pseudoinverse and the transpose of $C$, respectively. We index vectors and matrices as done in MATLAB; thus, the $k$ columns of $A$ with corresponding indices in vector $\bp \in \N_+^k$ are denoted by $A(:,\bp)$.

The DEIM algorithm \cite{Barrault,Chaturantabut} is a technique used to select some important column and row indices from an $m \times n$ data matrix $A$ \cite{Sorensen}, where without loss of generality  $m \ge n$. We wish to select $k\ll n$ relevant row and column indices. The first step in the DEIM procedure is to compute a (reduced) SVD, $A = U \Sigma V^\top$. The associated cost when a direct method is used is $\calo(mn^2)$, independent of the value of $k$. The paper also considers iterative methods to approximate the leading $k$ singular vectors. Having the left and right singular vectors contained in $U$ and $V$, respectively, we are interested in selecting distinct row indices $s_1, \dots, s_k$ from the set $\{1,\dots,m\}$ and column indices $p_1, \dots, p_k$ from the set $\{1,\dots,n\}$. The result of the method may also be represented by an $m \times k$ row {\em selection matrix} $S$ and an $n \times k$ column {\em selection matrix} $P$, whose columns are the standard basis vectors indexed by the selected indices.
The corresponding CUR factorization is (instead of the conventional use of the letter $U$ for the middle matrix we will use $M$ because $U$ is used to denote the matrix containing left singular vectors)
\begin{equation}\label{eq: cur}
\begin{array}{ccccc}
A& \approx & C & M & R~, \\[-0.5mm]
{\scriptstyle m\times n} & & {\scriptstyle m\times k} &{\scriptstyle k\times k} & {\scriptstyle k\times n}~
\end{array}
\end{equation} 
where the full-rank matrices $C=AP$ and $R=S^\top A$ consist of a subset of the columns and rows of A, respectively, and the middle matrix $M$ of full rank is computed such that the decomposition is as close to $A$ as possible. { For a general overview on computing the middle matrix, see \cite{hamm2021perturbations,hamm2020perspectives}}. Given $C$ and $R$, a standard procedure to determine $M$
(see, e.g., \cite[Sec.~2]{Sorensen}, where also an alternative is presented)
is by two consecutive least squares problems:

\begin{tabular}{ll}
{\footnotesize 1:} &
Solve the least squares problem $CX \approx A$ for $X \in \R^{k \times n}$ \\
& with solution $X = (C^\top C)^{-1} C^\top A$. \\
{\footnotesize 2:} &
Solve the least squares problem $R^\top M^\top  \approx X^\top$ for $M \in \R^{k \times k}$\\
& with solution $M = X R^\top  (RR^\top)^{-1}$.
\end{tabular}

Both steps are optimal with respect to the spectral and Frobenius norm. It is important to note that the solution in the spectral norm may not be unique. In many applications, one cares primarily about key columns or rows of $A$, rather
than an explicit $A\approx CMR$ factorization. Thus, an interpolative decomposition of the form $A = C\wh M$ or $A = \wh M R$ \cite{Voronin}.

We will now describe the process of selecting row indices using the DEIM scheme, and the procedure for choosing column indices is analogous. The DEIM algorithm begins with the leading left singular vector $\bu_1$, and the initial index, denoted as $s_1$, corresponds to the entry with the largest magnitude in $\bu_1$, i.e., $|\bu_1(s_1)|=\norm{\bu_1}_\infty$, where $\norm{\cdot}_\infty$ denotes the infinity-norm. Given that $I$ is the identity matrix, let $\bs=[s_1]$, $S_1=I(:,s_1)$, $U_1=[\bu_1]$, and define an oblique projection operator as $\mathbb{S}=\bu_1(S_1^\top \bu_1)^{-1}S_1^\top$. 

Suppose we have $j-1$ indices, so that 
\[\bs_{j-1}=\begin{bmatrix}s_1\\\vdots\\s_{j-1}\end{bmatrix}, \quad S_{j-1}=I(:,\bs_{j-1}), \quad U_{j-1}=[\bu_1,\dots,\bu_{j-1}],\]
and 
\[\mathbb{S}_{j-1}=U_{j-1}(S_{j-1}^\top U_{j-1})^{-1}S_{j-1}^\top.\]
Compute the residual vector $\br_j=\bu_j-\mathbb{S}_{j-1}\bu_j$, and choose the subsequent index $s_j$ such that $|\br_j(s_j)|=\norm{\br_j}_\infty$. It is important to point out that employing the oblique projection operator $\mathbb{S}_{j-1}$ on $\bu_j$ ensures that the $s_{j-1}$ entry in $\br_j$ remains 0, thereby guaranteeing non-repeating indices. Additionally, it is worth mentioning that, for the projector $\mathbb{S}_{j-1}$ to exist at the $j$th step, $S_{j-1}^\top U_{j-1}$ must be nonsingular, which is guaranteed by the linear independence of the columns in matrix U.
\cref{algo:DEIM} summarizes the DEIM index selection procedure \footnote{The backslash operator utilized in these algorithms follows a MATLAB-like convention for solving linear systems and least-squares problems.}. A variant of the DEIM scheme proposed by \cite{Drmac} called QDEIM involves computing a column-pivoted QR factorization on the transpose of the singular vectors to obtain the column and row indices, which correspond to the indices of the first $k$ pivoted columns.
\begin{algorithm}[htb!]
\KwData{$U \in \R^{m \times k}$ with $k\le m$ (of full rank)}
\KwResult{Indices $\bs \in \N_+^k$ with non-repeating entries}

$\bs(1)$ = $\argmax_{1\le i\le m}~ |(U(:,\,1))_i|$

\For{ $j = 2, \dots, k$}{
 $U(:,\,j) = U(:,\,j)-U(:,\,1:j-1)\cdot (U(\bs,\,1:j-1)\ \backslash \ U(\bs,\,j))$ 
 
 $\bs(j)$ = $\argmax_{1\le i\le m}~ |(U(:,\,j))_i|$\hspace{3mm}
}
\caption{DEIM index selection scheme \cite{Sorensen}}\label{algo:DEIM}
\end{algorithm}

In this paper, we explore iterative subselection variants of the DEIM scheme. Our contribution aims to improve the approximation quality of the DEIM scheme by iteratively invoking it, in the sense that we select subsequent columns and rows based on the previously selected ones. Thus, we modify $A$ after each iteration by removing the information that has been captured by the previously selected columns and rows. We show this by adapting an existing volume sampling technique for the DEIM scheme and also propose a new strategy. We also discuss how iterative procedures for computing a few singular vectors of large data matrices can be used with our proposed methods. To the best of our knowledge, this is the first deterministic DEIM type algorithm for large-scale data sets. 

\section{Volume sampling for column subset selection problem}\label{sec:adapsam}
The iterative subselection strategies proposed in this work are related to the so-called \emph{volume sampling} for column subset selection. In this section, we provide an overview of the volume sampling technique proposed by  \citeauthor{deshpande2006matrix} \cite{deshpande2006matrix}. The authors introduce a probabilistic method that iteratively selects a subset of columns in multiple rounds to construct a rank-$k$ approximation of a matrix. This approach has been demonstrated to provide improved accuracy and flexibility compared to \emph{one-round sampling} methods. One-round sampling methods refer to selection schemes that obtain all $k$ columns in a single round. 

\begin{algorithm}[htb!]
    \caption{Volume sampling for column subset selection \cite{deshpande2006matrix}}\label{algo:AdapSamp1}
    \KwData{$A\in\R^{m\times n}$, target rank $k$, \# rounds of $t$, columns per round $c$}
    
    \KwResult{$C\in\R^{m\times tc}$}

    $\bp=[\ ]; \quad E_0=A$

    \For{$i=1, \dots, t$}{

 \For {$j = 1, \dots, n$}{
 
     \textbf{if} $j\in\bp$ \ \textbf{then} \ $\text{pr}_j^{(i)}=0$\label{probline1} \hfill (sample without replacement)

     \textbf{else} $\text{pr}_i^{(j)}={\norm{E_{i-1}^{(j)}}^2}\,/\,{\norm{E_{i-1}}_F^2}$\label{probline}
    
}

$\bp_i$\, =\, set of $c$ indices sampled according to $\text{pr}_{i}$

$\bp= [\bp \ \ \bp_i]$

 $C=A(\,:,\bp)$; \quad $E_i=A-CC^+A$
    }
    
\end{algorithm}

The volume sampling method of \cite{deshpande2006matrix} as summarized in \cref{algo:AdapSamp1} involves alternating between two steps in each round: selecting a subset of columns and updating the probability distribution over all columns. The selection of columns in each round is influenced by the columns picked in previous rounds.  Suppose we aim to select a subset of $k$ columns from matrix $A$, the process begins with an initial probability distribution and randomly selects $c<k$ columns to form a matrix $C$. The selection of columns is based on the norms of the columns, as described in \cite{deshpande2006matrix,frieze2004fast}. Each column $j$ is chosen with a probability $\text{pr}_i^{(j)}={\norm{E_{i-1}^{(j)}}^2}\,/\,{\norm{E_{i-1}}_F^2}$ (as in \cref{probline} of \cref{algo:AdapSamp1}). After selecting $c$ columns, the probabilities are updated based on the chosen columns, and $c$ new columns are sampled and added to the matrix $C$. This iterative process continues until all $k$ columns are selected. Note that assigning zero probability to previously selected indices, as described in \cref{probline1}, represents a sampling without replacement strategy.

The authors present a detailed explanation and theoretical analysis of this volume sampling technique, emphasizing its advantages and diverse applications \cite{deshpande2006matrix}. The algorithm improves the accuracy of a CUR decomposition compared to one-round sampling methods as demonstrated in \cite{deshpande2006matrix,deshpande2006adaptive,paul2015column,Zhang}. Moreover, it allows for flexibility by accommodating different criteria for selecting column and row subsets based on specific problem requirements, which we will discuss in \cref{sec:iterDEIM}. In their approach (\cref{algo:AdapSamp1}), a constant number of columns is selected per iteration, and the residual is computed as $E=A-CC^+A$.

\pg{In this paper, we introduce a modified approach based on the DEIM scheme to implement a variant of the volume sampling proposed by \citeauthor{deshpande2006matrix} \cite{deshpande2006matrix} and propose a new iterative subselection strategy. By incorporating the DEIM scheme, our methods provide a deterministic technique for iterative subselection of column indices, in contrast to the original methods that employ a probabilistic approach \cite{optimalboutsidis,deshpande2006matrix,deshpande2006adaptive,paul2015column,Zhang}. 

We propose deterministic variants of the volume sampling technique for several reasons:
\begin{enumerate}
    \item The volume sampling technique, akin to many randomized algorithms requires oversampling of columns and rows beyond the specified target rank $k$ to achieve strong provable approximation guarantees. In our experience, the deterministic algorithms (including our proposed algorithms) typically yield lower approximation errors compared to randomized algorithms when we fix a rank parameter $k$ and choose exactly $k$ columns and rows (see, e.g., \cref{tab:exprand}). 
    \item For large-scale data sets, the proposed algorithms eliminate the need for explicitly computing the residual matrix $E=A-CC^+A$ as done in \cite{deshpande2006matrix}, making it efficient and suitable for large-scale matrices (see, e.g., \cref{tab:exprand}).
    \item Reproducibility is ensured with deterministic algorithms, unlike randomized sampling methods where results may not always be reproducible, even with a set seed. Note that, despite the use of an initial random vector in the Krylov--Schur routine of our large-scale deterministic algorithm, its influence is minimal, resulting in consistent and easily reproducible SVD.
\end{enumerate}}

\section{Small-scale DEIM type CUR with iterative SVDs}\label{sec:iterDEIM}
In this section, we introduce new index-picking schemes for constructing a CUR decomposition.
The standard DEIM scheme computes an SVD of $A$ once, after which the indices
are picked iteratively ``locally optimal''.
The new methods that we present now compute an SVD in every iteration. The algorithms adaptively select columns and rows of $A$ in several rounds. In each iteration, we modify $A$ by removing the information that has been captured by the previously selected columns and rows. The time complexities of the various proposed methods after $t$ rounds are summarized in \cref{ch6tab:overview}. This includes the computational time for computing an SVD and an updated $A$ (residual) matrix in every round. 

\begin{table}[htb!]
\centering \footnotesize
\caption{Summary of abbreviations for the various algorithms.}\label{tab:abbr}
{\begin{tabular}{l|lll}\hline\rule{0pt}{2.3ex}%
  Abbreviation& \# Indices per round & Residual&Algorithm\\ \hline\rule{0pt}{2.5ex}%
     \textbf{CADP-CX}&  Fixed number  & $A-CX$ &\cref{algo:AdapSamp}\\ 
     \textbf{CADP-CUR} & Fixed number  &  $A-CMR$& \cref{rmk2}\\
     \textbf{DADP-CX}&Singular value decay-based  & $A-CX$ &\cref{algo:ExAdapDEIM}\\
    \textbf{DADP-CUR}& Singular value decay-based & $A-CMR$ & \cref{algo:NewAdapDEIM}\\ \hline
    
\end{tabular}}
\end{table}

\begin{table}[htb!]
\centering
\caption{Summary of the dominant work of the different algorithms after $t$ rounds. The time complexity column excludes the computational cost of the DEIM scheme as it is approximately the same for all algorithms.}\label{ch6tab:overview}
{\scriptsize\begin{tabular}{lll|lll}\hline\rule{0pt}{2.3ex}%
Method                   & {Matrix} & {SVD}  & \multicolumn{3}{c}{Time }             \\ 
                         &      &                    & {svd}             & {$X$ or $M$} &  {Residual ($E$)}                        \\  \hline\rule{0pt}{2.5ex}%
\textbf{CADP-CX}      & \multirow{4}{*}{Small}        & \multirow{4}{*}{Full}        & \multirow{4}{*}{$\calo(tmn^2)$}&\multirow{4}{*}{$\calo(tmnk)$}   &\multirow{4}{*}{$\calo(tmnk)$}             \\
\textbf{DADP-CX}    &    &            & &   &              \\
\textbf{CADP-CUR}   &      &            &  &   &               \\
\textbf{DADP-CUR}   &      &            &  &   &               \\[0.5mm] \hdashline \rule{0pt}{3ex}%
\textbf{Large-scale: DADP-CX} & Large   & Few            & {$\calo(mn\cdot \text{nr}_\text{in})$} & $\calo(mnk)$    &\multicolumn{1}{c}{--} \\ \hline
\end{tabular}}
\end{table}

\pg{\begin{remark} \rm
The complexity estimates and the other descriptions in \cref{ch6tab:overview} are similar for the first four algorithms. However, when constructing a CUR factorization using the first two, one needs to compute almost twice the number of SVDs, $X$, and $E$, compared to what is required by \textbf{CADP-CUR} and \textbf{DADP-CUR}. For the large-scale algorithm, estimating the precise time complexity of computing a low-rank SVD using an iterative method can be challenging since it depends on the total number of inner iterations (denoted by $ \text{nr}_\text{in}$ in the table) needed.
\end{remark}}

{\subsection{An iterative DEIM with fixed indices per round and one-sided projected residual}}\label{sec:AdapSamp}
We present a deterministic variant of the iterative subselection scheme discussed in \cref{sec:adapsam}. The proposed algorithm called \textbf{CADP-CX} builds upon the original volume sampling algorithm \cite{deshpande2006matrix} by leveraging the benefits of the DEIM technique.  The method involves iteratively selecting a constant number of column indices, denoted as $c$, from $A$ in multiple rounds. We start by computing the leading $c<k$ singular vectors of $A$. Next, we apply the DEIM scheme (\cref{algo:DEIM}) to these singular vectors, resulting in the first set of $c$ indices. We then update $A$ by computing the residual matrix $E$ using the interpolative decomposition (as described in \cref{line:res} of \cref{algo:AdapSamp}). Next, we compute the leading $c$ singular vectors of $E$ and apply the DEIM procedure again to obtain the next set of $c$ indices. This process is repeated until we have selected all $k$ required indices. The procedure is summarized in \cref{algo:AdapSamp}.

\begin{algorithm}[htb!]
    \caption{Iterative DEIM with fixed indices per round and one-sided projected residual}\label{algo:AdapSamp}
    \KwData{$A\in\R^{m\times n}$, target rank $k$, columns per round $c$ (with $c\,|\,k$)}
    
    \KwResult{Low-rank CUR decomposition
$A_k = A(:,\bp) \, \cdot \, M \, \cdot \, A(\bs,:)$}

    Set $E=A$;\quad $\bp=[\ ]$;\quad $\bs=[\ ]$

    \For{$i=1, \dots, k/c$}{

     Compute $[\sim, \sim, V]$ = {\sf svd}($E$)

 $\bp_i={\sf deim}(V(:,1:c))$ \hfill (Iteratively pick $c$ indices) 
    
     $\bp= [\bp \ \ \bp_i]$

     $C=A(\,:,\bp)$; \quad $X=C\backslash A$;\quad  $E=A-CX$ \label{line:res}    
    }

    Repeat steps 1--7 on $A^\top$ to find the row indices $\bs$\label{repeat:steps1}

$M = A(:,\bp) \ \backslash \ (A \ / \ A(\bs,:))$
    
\end{algorithm}

With regards to the memory and computational complexity, computing the residual in \cref{{line:res}} involves a full iteration over the matrix, which has a space complexity of $\calo(mn)$. Given that we use the DEIM procedure and select $c$ columns per iteration, in terms of computational complexity, a full SVD requires $\calo(mn^2)$, one run of the DEIM algorithm requires $\calo(mc^2)$, and computing the residual in \cref{line:res} costs $\calo(mnk)$. The overall time complexity after $t$ rounds is $\calo(tmn^2 +tmc^2 +tmnk)$.

\subsection{A new iterative subselection method}\label{sec:ExAdapDEIM}
In \cref{algo:ExAdapDEIM}, we introduce a new iterative subselection strategy referred to as \textbf{DADP-CX} for a CUR factorization, which differs from the method employed in our proposed \cref{algo:AdapSamp} and the adaptive sampling procedures in previous works \cite{optimalboutsidis,deshpande2006matrix,paul2015column,Zhang}. In contrast to the previous strategy, which selects a fixed number of columns or rows in each iteration, this new strategy dynamically adjusts the selection schedule based on the decay of the singular values of the data (the relative magnitudes of the singular values).

The motivation behind this new approach is to adapt the subselection process according to the significance of the singular values. By considering the decay pattern of the singular values, we can prioritize the selection of columns or rows that contribute the most to the data's overall structure and information.
The decay of singular values provides valuable information about the significance of different components in the data. By leveraging this information, the iterative subselection strategy can adapt to the specific characteristics of the data and prioritize the selection of columns or rows that contribute the most to its structure. This adaptability allows for a more data-driven selection process.

\begin{algorithm}[htb!]
\caption{Singular value decay-based iterative DEIM with one-sided projected residual}\label{algo:ExAdapDEIM}

\KwData {$A \in \R^{m \times n}$, desired rank $k$, threshold parameter $\delta\in (0,1]$, upper limit $\ell$ }

\KwResult {Low-rank CUR decomposition
$A_k = A(:,\bp) \, \cdot \, M \, \cdot \, A(\bs,:)$ }
Set $E=A$;\quad $\bp=[\ ]$;\quad $\bs=[\ ]$

\While{{\sf length}$(\bp) <k$}{
 
 Compute $[\sim, \Sigma, V]$ = {\sf svd}($E$)\label{svdline1}
 
Let $b$ be the last index $i \le k-{\sf length}(\bp)$ with $\sigma_i \ge \delta \, \sigma_1$

$c=\min(b,\ell)$; \quad $\bp_c={\sf deim}(V(:,1:c))$  
 
 $\bp=[\bp \ \ \bp_c]$;\quad $C=A(:,\bp)$; \quad $X=C\backslash A$
 
 $E = A - CX$\label{residual}\hfill (Update matrix) 
}

Repeat steps 1--8 on $A^\top$ to find the row indices $\bs$\label{repeat:steps}

$M = A(:,\bp) \ \backslash \ (A \ / \ A(\bs,:))$

\end{algorithm}

With a user-defined threshold $\delta\in (0,1]$, the small-scale version of our method begins by computing the leading singular vectors corresponding to the singular values greater or equal to the threshold multiplied by the largest singular value of $A$, i.e., all $\sigma_i \ge \delta\cdot\sigma_1$.  Let $b$ denote the number of singular values satisfying this condition. Additionally, we introduce an extra parameter $\ell$ to establish an upper limit on the number of indices per round, taking into account the number of singular values that exceed the threshold. Consequently, we select the first $c=\min(b,\ell)$ column indices denoted by $\bp_c$ by applying the DEIM scheme to the leading $c$ right singular vectors ($V_c$).

We then proceed to construct an interpolative decomposition using the chosen column indices and compute the residual matrix $E$ by subtracting this approximation from $A$, i.e., $E=A-CC^+A$. To determine the next set of indices, we repeat the aforementioned process on $E$. Thus, we compute the leading singular vectors of $E$ corresponding to singular values greater than $\delta$ times the largest singular value of $E$ and repeat the procedure mentioned earlier. 
We expect that the multiple passes through $A$ would lead to a reduced approximation error. It is worth mentioning that in \cref{repeat:steps}, there is no need to compute the initial SVD of $A^\top$ since we can store the initial left singular vectors from the SVD of $A$.

In addition to the new selection strategy described in \cref{algo:ExAdapDEIM}, we also define an alternative way to compute the residual in the index selection process, which is presented in \cref{algo:NewAdapDEIM}. The newly proposed iterative subselection algorithms (\cref{algo:AdapSamp,algo:ExAdapDEIM}) and existing adaptive sampling procedures such as those outlined in \cite{optimalboutsidis,deshpande2006matrix,deshpande2006adaptive,paul2015column, Zhang} define the residual as the error incurred by projecting the matrix $A$ onto either the column space of $C$ or the row space of $R$, i.e., $E = A - CC^+A$ or $E = A - AR^+R$, respectively. In contrast, this new method termed as \textbf{DADP-CUR}, defines the residual as the error incurred by simultaneously projecting $A$ onto both the column space of $C$ and the row space of $R$. This means computing a CUR factorization at each step using only the selected columns and rows. 

{Note that, with this residual $E=A-CMR$, every column and row of  $A$ has a chance of being selected in subsequent rounds, irrespective of whether they were selected previously. Consequently, this might lead to the re-selection of columns and rows that were chosen in earlier rounds. \Cref{line:withoutrep} of \cref{algo:NewAdapDEIM} is a possible sample without-replacement strategy that can alleviate this problem. Since the DEIM procedure selects the index corresponding to the entry of the largest magnitude in a given vector, when these indices are zeroed out after being chosen, it guarantees that they will not be selected.}

\begin{algorithm}[htb!]
\KwData{ $A \in \R^{m \times n}$, desired rank $k$, threshold parameter $\delta=(0, 1]$, upper limit $\ell$}
\KwResult{ Low-rank CUR decomposition
$A_k = A(:,\bp) \, \cdot \, M \, \cdot \, A(\bs,:)$ }

Set $E=A$;\quad $\bp=[\ ]$;\quad $\bs=[\ ]$

\While{{\sf length}$(\bp) <k$}{

 $[U, \Sigma, V]$ = {\sf svd}($E$)\label{svdline2}
 
 Set $V(\bp,:\ )=0$, \quad $U(\bs,:\ )=0$ \label{line:withoutrep} 
 
 Let $b$ be the last index $i \le k-{\sf length}(\bp)$ with $\sigma_i > \delta \, \sigma_1$

$c=\min(b,\ell)$; \quad $\bp_c={\sf deim}(V(:,1:c))$  
 
 $\bs_c={\sf deim}(U(:,1:c))$

 $\bp=[\bp \ \ \bp_c]$, \quad $\bs=[\bs \ \ \bs_c]$
 
 $M = A(:,\bp) \ \backslash \ (A \ / \ A(\bs,:))$
 
 $E = A - A(:,\bp)\cdot M \cdot A(\bs,:) $ \label{residual2}\hfill (Update matrix) 
 
}
 \caption{Singular value decay-based iterative DEIM with two-sided projected residual}\label{algo:NewAdapDEIM}
\end{algorithm}

By considering the simultaneous projection onto the column and row spaces, we aim to use a residual that provides a more accurate representation of the error in the CUR factorization. It takes into account the combined effect of selecting specific columns and rows on capturing the underlying structure and information in the data. This approach offers several potential advantages. It allows for a more comprehensive assessment of the error in the index selection process, considering the contributions from both the columns and rows. Furthermore, it ensures that the residual accurately reflects the approximation quality obtained by a CUR factorization using the selected columns and rows. Additionally, it has the potential to reduce computational costs compared to  \cref{algo:ExAdapDEIM}, as the latter approach involves performing nearly twice the number of SVDs required by \cref{algo:NewAdapDEIM}. 

{\begin{remark}\label{rmk2} {\rm
For the completeness of comparison, in the experiments, we consider a variation of \cref{algo:AdapSamp} referred to as {\bf CADP-CUR}, where we use the newly defined residual. Thus, in \cref{algo:AdapSamp}, {\bf CADP-CUR} calculates the right and left singular values, selects both column and row indices in each round, and updates the residual using $E=A-CMR$, where $M = C^+AR^+$.}
\end{remark}}

Note that when $\delta=0$, both \cref{algo:ExAdapDEIM,algo:NewAdapDEIM} are equivalent to the DEIM type CUR factorization. In terms of time complexity, suppose we need $t$ iterations in \cref{algo:ExAdapDEIM,algo:NewAdapDEIM} to select all $k$ columns and rows. The cost of solving $E$ is $\calo(tmnk)$. The cost of an SVD and one run of the DEIM scheme are $\calo(mn^2)$ and $\calo(nc^2)$, respectively, where $c$ is the maximum number of columns selected per iteration.  Therefore, the overall cost of the algorithms is $\calo(tmn^2 +t(m+n)c^2 +tmnk)$. However, constructing $C$ and $R$ using \cref{algo:ExAdapDEIM} requires two runs of it. Thus, its cost is almost twice that of \cref{algo:NewAdapDEIM}.

For small matrices, these iterative subselection techniques may be worthwhile as the costs are modest and the quality of the approximations may increase. However, these schemes may be especially interesting for large matrices, for which an SVD may be too expensive, and iterative methods are used to compute the left and right singular vectors. We will study this situation in the next section.

\section{Large-scale DEIM type CUR with iterative SVDs}\label{sec:large-scale}
For large-scale matrices, taking an SVD every round in \cref{algo:AdapSamp,algo:ExAdapDEIM,algo:NewAdapDEIM} will usually be prohibitively expensive. Indeed, even one (reduced) SVD will be too costly, which means that
the standard DEIM type CUR decomposition is generally not affordable.
However, the proposed algorithm is suitable for large-scale data, as approximating the largest singular vectors by iterative (Krylov) methods is usually a relatively easy task. Additionally, here, we do not explicitly compute the residual matrix as done in the proposed algorithms; this is done implicitly in the computation of the approximate singular vectors. Furthermore, instead of computing the full SVD as we do in \cref{algo:AdapSamp,algo:ExAdapDEIM,algo:NewAdapDEIM}, we now carry out: \\[1mm]
\begin{tabular}{ll}
{\footnotesize 1:} & \phantom{M} Approximate $\wh U$ and $\wh V$ of $E$.
\end{tabular}

{This can efficiently be carried out by Krylov--Schur for the SVD \cite{stoll2012krylov}, a very efficient implicitly restarted version of Lanczos bidiagonalization (which in our experience is generally considerably faster than the mathematically equivalent method of \cite{baglama2005augmented} as implemented in Matlab's \texttt{svds})}. The idea is as follows. Let $k < {\wh k}$ be the minimal and maximal dimension of the subspaces. We first carry out ${\wh k}$ steps of Lanczos bidiagonalization summarized by the matrix equations
\[
E\wh \, V_{\wh k} = \wh U_{\wh k} \, B_{\wh k}, \quad
E^\top \wh U_{\wh k}= \wh V_{\wh k}\, B_{\wh k}^\top  + \beta_{\wh k}\, \wh \bv_{{\wh k}+1} \, \be_{\wh k}^\top,
\]
where $B_{\wh k}$ is bidiagonal. The singular values of $B_{\wh k}$ are approximations to those of $E$,
and the singular vectors lead to approximations to those of $E$.
With the SVD $B_{\wh k}= W \wh \Sigma Z^\top$, we get
\[
E\, (\wh V_{\wh k}Z) = (\wh U_{\wh k}W) \, \wh \Sigma, \quad
E^\top(\wh U_{\wh k}W) = (\wh V_{\wh k}Z) \, \wh \Sigma + \beta_{\wh k}\wh \bv_{{\wh k}+1} (W^\top  \be_{\wh k})^\top.
\]
For any upper triangular matrix $\wh \Sigma$ an elegant implicit restart procedure
is possible; here $\wh \Sigma$ is even diagonal.
Order the singular values in the desired way; in this case nonincreasingly.
Partition the transformed basis, redefining $\wh U_k$ and $\wh V_k$:
\begin{equation}
\label{partition}
\wh U_{\wh k}W =: [\wh U_k \ \, \wh U_{{\wh k}-k}], \quad
\wh V_{\wh k}Z =: [\wh V_k \ \, \wh V_{{\wh k}-k}], \quad
\wh \Sigma = \smtxa{cc}{\wh \Sigma_k \\ & \wh \Sigma_{{\wh k}-k}},
\end{equation}
and redefine $B_k = \wh \Sigma_k$, $\beta_{k+1} := \beta_{{\wh k}+1}$,
$\wh \bv_{k+1} := \wh \bv_{{\wh k}+1}$, and $\bff_k := W^\top  \be_{\wh k}$.
We can now conveniently restart from the decomposition
\[
E\wh V_k = \wh U_k \, B_k, \quad
E^\top  \wh U_k = \wh V_k B_k^\top  + \beta_k \, \wh\bv_{k+1} \, \bff_k^\top.
\]
The pair $(\wh U_k, \wh V_k)$ may be viewed as a pair of approximate invariant spaces with error $\|\bff_k\|$. The spaces are expanded with Lanczos bidiagonalization to dimension ${\wh k}$, after which the selection procedure is carried out again. This scheme is repeated until the quantity $\|\bff_k\|$ is sufficiently small. We summarize the method in \cref{algo:itersvd}. Note that MATLAB built-in function {\sf svds} is a different implementation of a related technique.

\begin{algorithm}[htb!]
\KwData{$E \in \R^{m \times n}$, desired rank $k$, initial vector $\bv_1$,
minimum and maximum dimension $k < {\wh k}$, tolerance {\sf tol}}

\KwResult{Approximation to $k$ largest singular triplets
$(\sigma_i, \bu_i, \bv_i)$, giving best low-rank approximation
$E_k = \wh U_k \wh \Sigma_k \wh V_k^\top$}

Generate $E\wh V_k = \wh U_k B_k$, \ $E^\top \wh U_k = \wh V_kB_k^\top  + \beta_{k+1} \wh \bv_{k+1} \bff_k^\top$

\For{$i = 1, 2, \dots$}{
 Expand to $E\wh V_{\wh k}= \wh U_{\wh k}B_{\wh k}$, \ $E^\top \wh U_{\wh k}= \wh V_{\wh k}B_{\wh k}^\top  + \beta_{{\wh k}+1} \wh \bv_{{\wh k}+1} \bff_{\wh k}^\top$ 

 Determine SVD $B_{\wh k}= W \wh \Sigma Z^\top$

 Partition according to \eqref{partition}, restart with $\wh U_k$, $\wh V_k$,

 \quad redefining $B_k := \wh \Sigma_k$, $\beta_{k+1} := \beta_{{\wh k}+1}$, $\wh \bv_{k+1} := \wh \bv_{{\wh k}+1}$, $\bff_k := W^\top  \be_{\wh k}$ 

 Stop if $\|\bff_k\| \le $ {\sf tol}
}

\caption{Krylov--Schur for the SVD \cite{stoll2012krylov}.}\label{algo:itersvd}
\end{algorithm}

\begin{algorithm}[htb!]
\caption{Large-scale: Singular value decay-based iterative DEIM with one-sided projected residual}\label{algo:NewAdapDEIM2}
\KwData {$A \in \R^{m \times n}$, desired rank $k$, threshold parameter $\delta \in (0,1]$, upper limit $\ell$}
\KwResult {Low-rank CUR decomposition
$A_k = A(:,\bp) \, \cdot \, M \, \cdot \, A(\bs,:)$ }
Set $E=A$;\quad $\bp=[\ ]$; \\
\While{{\sf length}$(\bp) <k$}{
Compute $\sigma_j$'s and $\bv_j$'s by \cref{algo:itersvd} \label{linesvds} \\
\qquad finding $b$, the last index $i \le k-{\sf length}(\bp)$ with $\sigma_i > \delta \, \sigma_1$ or \\
\qquad at most $\sigma_{\ell}$. Let $c=\min(b,\ell)$\\
$V(\bp,:\ )=0$;\quad  $\bp_c={\sf deim}(V(:,1:c))$ \\
$\bp=[\bp \ \ \bp_c]$; \quad $C=A(:,\bp)$ \\
Update an incremental QR decomposition $C=QT$ \label{lineincQR}\\
{$E$ is the function $\by = E(\bx, \text{transp\_flag})$ with: \\
\quad $\text{transp\_flag} = \texttt{true} : \ \ \by = A\bx; \ \by = \by - Q(Q^\top \by);$ \hfill \\ \label{implicit_E}
\quad $\text{transp\_flag} = \texttt{false} : \ \by = \bx - Q(Q^\top \bx); \ \by = A^\top \by$} \\
}
Repeat steps 1--9 on $A^\top$ to find the row indices $\bs$ \\
$M = A(:,\bp) \ \backslash \ (A \ / \ A(\bs,:))$
\end{algorithm}

In \cref{algo:NewAdapDEIM2} we provide the large-scale version of \cref{algo:ExAdapDEIM} by employing \cref{algo:itersvd}.  Without showing details, it is worth noting that this can also be adapted for \cref{algo:AdapSamp,algo:NewAdapDEIM}. It is important to note that the threshold parameter $\delta$ and the upper limit $\ell$ on the number of indices to be selected per round are incorporated within the implementation of \cref{algo:itersvd}. 

The efficiency of \cref{algo:NewAdapDEIM2} is because for the procedure in \cref{algo:itersvd} we do not need the matrix $E$ in explicit form; only matrix-vector products (MVs) with $E$ and $E^\top$ are necessary (see \cref{implicit_E}). The routine of \cref{algo:itersvd} also takes several MVs that depend on the distribution of the singular value and the starting vector. The cost of computing the singular values and vectors in \cref{linesvds} depends on the total number of inner iterations of \cref{algo:itersvd}. The number of iterations required by the Krylov--Schur algorithm depends on the size of the matrix. In a matrix-vector product $E\bx$ for a vector $\bx$, the component $A\bx$ costs $\calo(mn)$ for a full matrix. In the case of a sparse matrix with $d$ nonzeros entries per row, the cost reduces to $\calo(md)$. The aggregated cost of \cref{lineincQR} is only $\calo(mk^2)$. In \cref {implicit_E}, the computation of $Q(Q^\top \by)$ requires $\calo(mk)$ operations. The cost of solving the least squares problem $M = C^+AR^+$ would be $\calo(mnk)$, which is relatively expensive. Nevertheless, it is important to highlight that this step is necessary for all CUR methods as the final step.

As previously mentioned, it is not necessary to compute the initial SVD of $A^\top$ in this case, as we can simply retain the initial left singular vectors obtained from the SVD of $A$. The value of $\delta$ will typically depend on the data set.
A value close to $1$ may be favorable for the approximation result but is more expensive since \cref{algo:itersvd} needs to be carried out approximately $k$ times. However, having $\delta=1$ implies that we select one index per iteration, and thus, we need just the first right and left singular vectors of $E$ corresponding to the largest singular value. We can reduce the computational cost by specifying an earlier convergence criterion for finding the approximate leading right and left singular vectors. We use Wedin's theorem for this. The theorem bounds the distance between subspaces and the proof is in (cf., e.g., \cite[pp.~260--262]{stewart1990matrix}).

\begin{theorem}\label{thm:wedins}(Wedin's Theorem)
    Given $E\in\R^{m\times n}$, let 
    \[ [U_1\ U_2 \ U_3]^\top \ E\ [V_1 \ V_2]=\mtxa{cc}{\Sigma_1&0\\0&\Sigma_2\\0&0},\]
    be the SVD of $E$ (where the singular values are not necessarily nonincreasing). The singular subspaces of interest are in the column spaces of $U_1$ and $V_1$. Let the inexact/approximate singular subspaces be in the column spaces of $\wh U_1$ and $\wh V_1$ in the decomposition
    \[ [ \wh U_1\ \wh U_2 \ \wh U_3]^\top \ \wh E\ [\wh V_1 \ \wh V_2]=\mtxa{cc}{\wh \Sigma_1&0\\0& \wh \Sigma_2\\0&0}.\]
    Now let $\Phi$ be the matrix of canonical angles between $\Range(U_1)$ and $\Range(\wh U_1)$, and $\Theta$ be the matrix of canonical angles between $\Range(V_1)$ and $\Range(\wh V_1)$. Given the residuals $F_1=E\wh V_1 -\wh U_1\wh \Sigma_1$, $F_2=E^\top \wh U_1 -\wh V_1\wh \Sigma_1$, suppose that there is a number $\alpha>0$ such that 
    \[\min |\sigma(\wh \Sigma_1)-\sigma(\Sigma_2)| \ge \alpha \quad \text{and} \quad \sigma_{\min}(\wh \Sigma_1)\ge \alpha.\]
    Then 
\[
\sqrt{\smash[b]{\norm{\sin{\Phi}}_F^2 + \norm{\sin{\Theta}}_F^2}} \le \alpha^{-1} \, \sqrt{\smash[b]{\norm{F_1}_F^2+\norm{F_2}_F^2}}.
 \]    
\end{theorem}
\cref{thm:wedins} shows that the computed singular vectors extracted by the projection method are optimal up to the factor in the right-hand side of the above inequality. This implies that any change in the entries of the computed singular vectors is bounded by this factor. Note that $\Sigma_2$ is unknown. For our context, we use $\wh \Sigma_2$ as an approximation to $\Sigma_2$. Since we are only concerned with approximating the first leading right singular vector $\wh \bv_1$, we approximate $\alpha\approx\wh \sigma_1 - \wh \sigma_2$. Let $m_1(\wh \bv_1)$ and $m_2(\wh \bv_1)$ denote the largest and second-largest entries in $\wh \bv_1$, respectively, and let $\mathbf{f}_2=E^\top \wh \bu_1-\wh \sigma_1 \wh \bv_1$ be the residual vector (associated with residual matrix $F_2$). The above iterative SVD routine results in $\mathbf{f}_1=E\wh \bv_1 -\wh \sigma_1 \wh \bu_1 =0$ (associated with residual matrix $F_1$).
The DEIM algorithm selects the index corresponding to the largest element in the magnitude of a vector. Therefore, when $\delta=1$, one can set an early convergence criterion to find the first singular vector that corresponds to the largest singular value, using the following approximate bound:
\[m_1(\wh \bv_1)-m_2(\wh \bv_1)\lesssim 2\ (\wh \sigma_1 - \wh \sigma_2)^{-1} \ \norm{\mathbf{f}_2}.\]

\section{Error Analysis} 
In this section, we present a well-established theoretical error bound that pertains to a broad category of CUR factorizations. Consequently, this bound remains valid for the methods we propose. We refer the reader to \cite[Section~4]{Sorensen} for  a comprehensive, constructive proof. 

Let $P\in \R^{n\times k}$ and $S\in \R^{m\times k}$ be selection matrices with some columns of the identity indexed by the indices chosen by employing the index selection techniques proposed in this paper. 
\begin{theorem}\label{pp1} \cite[Thm.~4.1]{Sorensen}
Given $A\in \R^{m \times n}$ and a target rank $k$, let $V\in \R^{n\times k}$ and $U\in \R^{m\times k}$ be the leading $k$ right and left singular vectors of $A$. Suppose $C=AP$ and $R=S^\top \!A$ are of full rank, and $V^\top \!P$ and $S^\top U$ are nonsingular, then, with $M=C^+\!AR^+$, a rank-$k$ CUR decomposition constructed by the proposed techniques satisfies
\[\norm{A-CMR}\le(\eta_\bs + \eta_\bp)\,\sigma_{k+1} \quad with \quad \eta_\bs < \sqrt{\tfrac{nk}{3}}\,2^k~, \quad \eta_\bp < \sqrt{\tfrac{mk}{3}}\,2^k,\]
where $\eta_\bp=\norm{(V^\top\!P)^{-1}}$, $\eta_\bs=\norm{(S^\top U)^{-1}}$.
\end{theorem}

\section{Experiments}
We conduct numerical experiments to evaluate the empirical performance of the DEIM scheme \cite{Sorensen}, the QDEIM procedure \cite{Drmac}, the {\sf MaxVol} method \cite{Goreinov2010}, and the iterative subselection techniques discussed in this paper. The following are the iterative subselection methods we evaluate:
\begin{enumerate}
    \item [] \textbf{CADP-CX}: \ refers to \cref{algo:AdapSamp}.
    \item [] \textbf{DADP-CX}: \ represents \cref{algo:ExAdapDEIM}.
    \item [] \textbf{DADP-CUR}: \ corresponds to \cref{algo:NewAdapDEIM}.
    \item [] \textbf{CADP-CUR}:\ denotes the adapted version of \cref{algo:AdapSamp} with the residual defined as $E=A-CMR$.

\end{enumerate}
 To assess the effectiveness of our algorithms, we test them on various data matrices from different application domains, such as economic modeling, categorizing text and retrieving information, and computational fluid dynamics. Our evaluation includes synthetic and real-world data matrices, both sparse and dense, with varying sizes ranging from small to large scale, which we summarize in \cref{tab:exps}. In the implementation, we use the in-built Matlab functions {\sf qr} for the column pivoted QR and  {\sf svds} or {\sf svd} for the SVD computation, where the latter is used for small-scale matrices. For \cref{algo:itersvd,algo:NewAdapDEIM2}, we use our implementation of the Krylov--Schur method of \cite{stoll2012krylov} by incorporating the threshold parameter $\delta$ and upper limit $\ell$ on the number of singular vectors to be computed. Unless otherwise stated, in all the experiments we use as default the number of rounds $t=10$, the parameter $\delta=0.8$, and upper limit $\ell=k/10$.
\begin{table}[htb!]
\centering
\caption{Various examples and dimensions considered.\label{tab:exps}}
{\footnotesize
\begin{tabular}{cllcccc} \hline \rule{0pt}{2.3ex}%
{\bf Exp.} & {\bf Domain} & {\bf Matrix} & $m$ & $n$ \\ \hline \rule{0pt}{2.3ex}%
1 & Synthetic & ~Sparse & 100000 &\phantom{12}300\\
2 & Text categorization& ~Sparse & \phantom{123}139 &15210 \\
3 & Text categorization& ~Sparse & \phantom{12}8293 &18933 \\
4 & Economic modeling & ~Sparse & \phantom{1}29610 &29610\\
5 & Computational fluid dynamics & ~Sparse & \phantom{1}30412 &30412\\ \hline
\end{tabular}}
\end{table}

\begin{experiment}{\rm
In our first experiment, we investigate how different choices of $\delta$ and the number of rounds $t$ affect the approximation accuracy of the various iterative subselection strategies. We use the relative approximation error ${\norm{A-CMR}}\,/\,{\norm{A}}$ as the evaluation metric.
For this experiment just as in \cite{Sorensen}, we generate a sparse, nonnegative matrix $A \in \R^{m \times n}$, with $m=100000$ and $n=300$, of the form 
\[A=\sum_{j=1}^{10}\frac{2}{j}\, \bx_j\ \by_j^\top  + \sum_{j=11}^{300}\frac{1}{j}\, \bx_j\ \by_j^\top,\]
where $\bx_j \in \R^{m}$ and $\by_j \in \R^{n}$ are sparse vectors with random nonnegative entries (i.e., $\bx_j={\sf sprand}(m,1,0.025)$ and $\by_j={\sf sprand}(n,1,0.025)$). 

From \cref{fig: 1} we observe that increasing the number of rounds $t$ or $\delta$ does not necessarily lead to a monotonic decrease in the approximation errors in the 2-norm.
The result implies that one needs to carefully choose the parameter $\delta$ or the number of rounds to get the optimum advantage of using the iterative subselection strategies. For this experiment, we also observe that using the residual $E=A-CMR$ instead of $A-CC^+A$ for the iterative subselection yields better approximation errors in the delta strategy while it produces worse approximation errors in the constant number of columns strategy. 

\begin{figure}[htb!]
\centering
{{
%
%
\definecolor{mycolor1}{rgb}{0.00000,0.44700,0.74100}%
\definecolor{mycolor2}{rgb}{0.85000,0.32500,0.09800}%
\begin{tikzpicture}

\begin{groupplot}[
group style={
    group name=my plots,
    group size=2 by 2,
    xlabels at=edge bottom,
    ylabels at=edge left,
    horizontal sep=1.3cm,vertical sep=2cm,
    },
    ylabel= {$\| A - CMR \|\ /\ \| A \|$},
width=0.5\linewidth
]

\nextgroupplot[xmin=2,xmax=10,ymin=0.0215,ymax=0.024,xlabel={\# of rounds}]
\addplot [color=mycolor1,line width=1.5pt, mark=asterisk, mark options={solid, mycolor1}]
  table[row sep=crcr]{%
2	0.023153439648361\\
3	0.023315181525076\\
5	0.0237494956627834\\
6	0.0222375395566711\\
10	0.0227545072848429\\
};
\addlegendentry{CADP-CUR}

\addplot [color=mycolor2,line width=1.5pt, mark=o, mark options={solid, mycolor2}]
  table[row sep=crcr]{%
2	0.0238409695422888\\
3	0.0222719677620885\\
5	0.0231506151028399\\
6	0.0217619854053936\\
10	0.022647291429044\\
};
\addlegendentry{CADP-CX}

\nextgroupplot[xmin=0.4,xmax=0.8,xtick={0.4,0.5,0.6,0.7,0.8},ymin=0.0215,ymax=0.023,xlabel={Deltas ($\delta$s)}]

\addplot [color=mycolor1,line width=1.5pt, mark=square, mark options={solid, mycolor1}]
  table[row sep=crcr]{%
0.4	0.0221887872729132\\
0.5	0.0217620350240561\\
0.6	0.0216661866902728\\
0.7	0.0219817500733099\\
0.8	0.022027578090881\\
};
\addlegendentry{DADP-CUR}

\addplot [color=mycolor2,line width=1.5pt, mark=triangle, mark options={solid, rotate=180, mycolor2}]
  table[row sep=crcr]{%
0.4	0.0223916729331237\\
0.5	0.0228535316073208\\
0.6	0.022320314049708\\
0.7	0.0220521711649708\\
0.8	0.0221189745235457\\
};
\addlegendentry{DADP-CX}

\end{groupplot}

\end{tikzpicture}%
}}

\caption{Relative approximation errors for the various iterative subselection DEIM CUR approximation algorithms 
for $k = 30$. The right figure represents selecting a constant number of columns and rows per iteration; the left is the delta strategy. In all cases, increasing the number of rounds or delta does not lead to a monotonic decrease in the approximation errors.
\label{fig: 1}}
\end{figure}
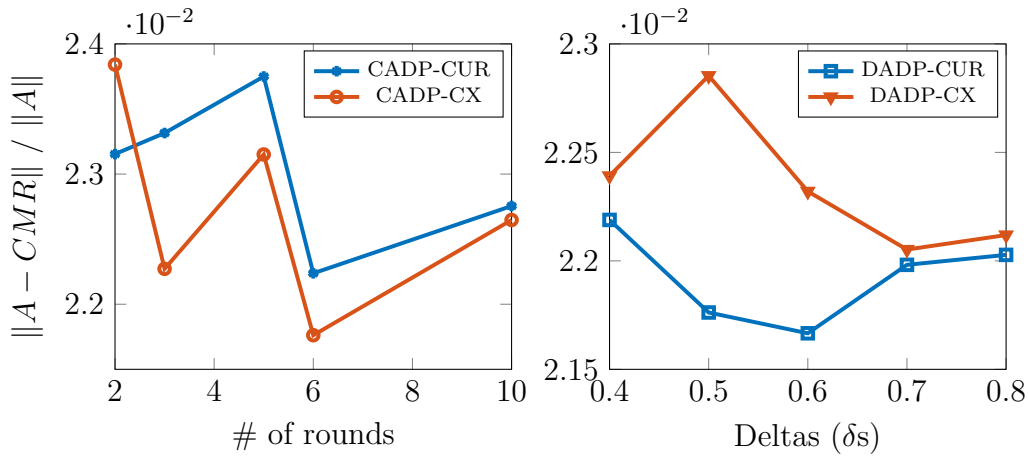
}
\end{experiment}

\begin{experiment}{\rm
Our next experiment is to demonstrate that the iterative subselection techniques yield better approximation results than one-round sampling. We perform the experiment using four real data sets and report the relative approximation error ${\norm{A-CMR}}\,/\,{\norm{A}}$ of each algorithm on each data set.

The first two data sets are relevant to text categorization and information retrieval applications. In such data analysis problems, a ``bag of words'' approach is commonly employed to represent documents. We opt for the {\sf Reuters-21578} text categorization collection, which comprises documents that were featured on Reuters' newswire in $1987$. This data set is extensively used as a benchmark in the text classification community, consisting of $21578$ documents categorized into $135$ categories. For our experiment, we use the preprocessed data set, which has $18933$ unique terms and $8293$ documents \cite{cai2005document}. We normalize the rows of the sparse matrix, which has dimensions $8293 \times 18933$, to have a unit length. The second data set, the Internet term document data, is from the Technion Repository of Text Categorization Datasets ({\sf TechTC}) \cite{Gabrilovich}. We use the test set $26$, which consists of a collection of $139$ documents on two topics (i.e., Evansville, Indiana (id 10567)
and Miami, Florida (id 11346)) with $15210$ terms describing each document\footnote{\url{http://gabrilovich.com/resources/data/techtc/}}. As in \cite{Sorensen}, the $139 \times 15210$ TechTC matrix rows are scaled to have a unit 2-norm.

The final two data sets are large sparse data matrices obtained from the publicly available SuiteSparse Matrix Collection.  One of these data sets, known as {\sf g7jac100}, results from the ``Overlapping Generations Model''  employed in the study of the social security systems of the G7 nations. This matrix is of dimensions $29610 \times 29610$, containing 335972 numerically nonzero entries and exhibiting a low rank of 21971. In addition, we use the matrix labeled {\sf invextr1-new}, which is linked to a computational fluid dynamics problem. The matrix has dimensions $30412\times 30412$ with a rank-$27502$ structure and contains $1793881$ nonzero entries.

From \cref{fig: 2}, we can see that in all cases our iterative subselection-based CUR algorithms have lower approximation errors than all the \emph{one-round} deterministic index selection algorithms considered. We also observe that the approximation error of the QDEIM and the {\sf MaxVol} techniques do not always decrease monotonically with increasing $k$ values. {This phenomenon is local to the spectral norm; from results not presented here the errors in Frobenius norm are monotonically non-increasing}. We choose the number of rounds for the CADP-CUR and CADP-CX algorithms to be $t=10$, and the parameter $\delta=0.8$ and upper limit $\ell=k/10$ for the DADP-CUR and DADP-CX algorithms. The results of all four proposed iterative subselection algorithms are comparable.   
\begin{figure}[htb!]
\centering
{{
%
%
\definecolor{mycolor1}{rgb}{0.00000,0.44700,0.74100}%
\definecolor{mycolor2}{rgb}{0,0,1}%
\definecolor{mycolor3}{rgb}{0.92900,0.69400,0.12500}%
\definecolor{mycolor4}{rgb}{0.49400,0.18400,0.55600}%
\definecolor{mycolor5}{rgb}{0.46600,0.67400,0.18800}%
\definecolor{mycolor6}{rgb}{0.30100,0.74500,0.93300}%
\definecolor{mycolor7}{rgb}{0,0,0}%
\begin{tikzpicture}

\begin{groupplot}[
group style={
    group name=my plots,
    group size=2 by 2,
    xlabels at=edge bottom,
    ylabels at=edge left,
    horizontal sep=1.3cm,vertical sep=2cm,
    },
    ymin=0.15,
    ymax=0.45,
    scaled y ticks=base 10:\exponent,
    xlabel=rank-$k$, 
    ylabel= {$\| A - CMR \|\ /\ \| A \|$},
legend style={at={(-0.2,-0.65)},anchor=south,legend  columns =4},
width=0.45\linewidth,
height=5cm
]
\nextgroupplot[title={\sf TechTC},xmin=10,xmax=50,xtick={10,20,30,40,50},ymin=0.35,ymax=0.7]
\addplot [color=red,line width=1.5pt, mark=diamond, mark options={solid, red}]
  table[row sep=crcr]{%
10	0.639521912990993\\
20	0.580782605356578\\
30	0.532279262004965\\
40	0.513507721559631\\
50	0.505725166948345\\
};

\addplot [color=mycolor6,line width=1.5pt, mark=triangle, mark options={solid, rotate=270, mycolor6}]
  table[row sep=crcr]{%
10	0.653243835239604\\
20	0.615415179609538\\
30	0.539511434041749\\
40	0.546096916853861\\
50	0.495024460797222\\
};
\addplot [color=mycolor7,line width=1.5pt, mark=triangle, mark options={solid, mycolor7}]
  table[row sep=crcr]{%
10	0.672793045850886\\
20	0.600829661691332\\
30	0.546136624140795\\
40	0.510254244300746\\
50	0.510247166844653\\
};
\addplot [color=mycolor2,line width=1.5pt, mark=asterisk, mark options={solid, mycolor2}]
  table[row sep=crcr]{%
10	0.572466808540772\\
20	0.489288177443187\\
30	0.424326383229177\\
40	0.402579590206831\\
50	0.376932450168752\\
};
\addplot [color=mycolor3,line width=1.5pt, mark=o, mark options={solid, mycolor3}]
  table[row sep=crcr]{%
10	0.572466808540772\\
20	0.471051461188483\\
30	0.427193142140243\\
40	0.390494073293096\\
50	0.368743503328529\\
};
\addplot [color=mycolor4,line width=1.5pt, mark=square, mark options={solid, mycolor4}]
  table[row sep=crcr]{%
10	0.570295922080816\\
20	0.475504213220872\\
30	0.440027481828186\\
40	0.397833578247394\\
50	0.374746744189983\\
};
\addplot [color=mycolor5, line width=1.5pt,mark=triangle, mark options={solid, rotate=180, mycolor5}]
  table[row sep=crcr]{%
10	0.570295922080816\\
20	0.475837525215249\\
30	0.426959321172227\\
40	0.390393695428244\\
50	0.366307310278865\\
};

\nextgroupplot[title={\sf Reuters},xmin=10,xmax=50,xtick={10,20,30,40,50},ymin=0.2,ymax=0.45]
\addplot [color=red, line width=1.5pt, mark=diamond, mark options={solid, red}]
  table[row sep=crcr]{%
10	0.419252145192696\\
20	0.40216393436931\\
30	0.350182743753939\\
40	0.303076873515701\\
50	0.299616838544632\\
};
\addplot [color=mycolor6, line width=1.5pt,  mark=triangle, mark options={solid, rotate=270, mycolor6}]
  table[row sep=crcr]{%
10	0.418502768333143\\
20	0.399715908766243\\
30	0.41341658280015\\
40	0.301295404182039\\
50	0.320429864424903\\
};
\addplot [color=mycolor7, line width=1.5pt, mark=triangle, mark options={solid, mycolor7}]
  table[row sep=crcr]{%
10	0.418253095775041\\
20	0.400016248445016\\
30	0.343211415455115\\
40	0.355826138602157\\
50	0.33240621125416\\
};

\addplot [color=mycolor2, line width=1.5pt, mark=asterisk, mark options={solid, mycolor2}]
  table[row sep=crcr]{%
10	0.364713755975854\\
20	0.298855176639389\\
30	0.257572601101225\\
40	0.230274878180724\\
50	0.216902387039248\\
};
\addplot [color=mycolor3, line width=1.5pt, mark=o, mark options={solid, mycolor3}]
  table[row sep=crcr]{%
10	0.364713755975854\\
20	0.293399313144267\\
30	0.254779434003222\\
40	0.230522076766072\\
50	0.212170266325728\\
};
\addplot [color=mycolor4, line width=1.5pt,  mark=square, mark options={solid, mycolor4}]
  table[row sep=crcr]{%
10	0.365029046095007\\
20	0.296033926805428\\
30	0.251682756505136\\
40	0.226640935736928\\
50	0.213593439530647\\
};
\addplot [color=mycolor5, line width=1.5pt,  mark=triangle, mark options={solid, rotate=180, mycolor5}]
  table[row sep=crcr]{%
10	0.365029046095007\\
20	0.287519438003208\\
30	0.246599940620695\\
40	0.222771894977483\\
50	0.205810610269483\\
};

\nextgroupplot[title={\sf g7jac100},ymin=0.05,ymax=0.55,xmin=100,xmax=500]
\addplot [color=red, line width=1.5pt,mark=diamond, mark options={solid, red}]
  table[row sep=crcr]{%
100	0.518471507862923\\
200	0.369512550201253\\
300	0.251480219806051\\
400	0.209715219208256\\
500	0.184295592438967\\
};

\addplot [color=mycolor6, line width=1.5pt,mark=triangle, mark options={solid, rotate=270, mycolor6}]
  table[row sep=crcr]{%
100	0.541811236388936\\
200	0.451919120462753\\
300	0.508499732344035\\
400	0.307062648095574\\
500	0.260784596094849\\
};

\addplot [color=mycolor7, mark=triangle, mark options={solid, mycolor7}]
  table[row sep=crcr]{%
100	0.549359961674411\\
200	0.451462467689926\\
300	0.508501801217514\\
400	0.290841520527453\\
500	0.253102550586343\\
};

\addplot [color=mycolor2,line width=1.5pt, mark=asterisk, mark options={solid, mycolor2}]
  table[row sep=crcr]{%
100	0.288166348705356\\
200	0.204583759719981\\
300	0.148593211509307\\
400	0.120796732558351\\
500	0.106659109546762\\
};

\addplot [color=mycolor3,line width=1.5pt, mark=o, mark options={solid, mycolor3}]
  table[row sep=crcr]{%
100	0.282223510130486\\
200	0.18573883880983\\
300	0.140477612418826\\
400	0.11460655536138\\
500	0.0986276276058334\\
};

\addplot [color=mycolor4,line width=1.5pt, mark=square, mark options={solid, mycolor4}]
  table[row sep=crcr]{%
100	0.293644461104449\\
200	0.20199534572069\\
300	0.15340342338346\\
400	0.123202018902749\\
500	0.105205070049218\\
};

\addplot [color=mycolor5,line width=1.5pt, mark=triangle, mark options={solid, rotate=180, mycolor5}]
  table[row sep=crcr]{%
100	0.28967994469243\\
200	0.192686554218367\\
300	0.138052882601562\\
400	0.114810278396209\\
500	0.0875765273236556\\
};

\nextgroupplot[title={\sf invextr1-new},ymin=0.05,ymax=0.16,xmin=500,xmax=1000,xtick={500,600,700,800,900,1000}]
\addplot [color=red,line width=1.5pt, mark=diamond, mark options={solid, red}]
  table[row sep=crcr]{%
500	0.151306884818792\\
600	0.14282243851067\\
700	0.120232246524414\\
800	0.109750582006787\\
900	0.0921435500111651\\
1000	0.0902377051844859\\
};
\addlegendentry{DEIM}

\addplot [color=mycolor6,line width=1.5pt, mark=triangle, mark options={solid, rotate=270, mycolor6}]
  table[row sep=crcr]{%
500	0.15240195676784\\
600	0.143920494187863\\
700	0.126076934877796\\
800	0.11010670687824\\
900	0.0975886625515651\\
1000	0.0934907666301426\\
};
\addlegendentry{MAXVOL}

\addplot [color=mycolor7,line width=1.5pt, mark=triangle, mark options={solid, mycolor7}]
  table[row sep=crcr]{%
500	0.152402473035168\\
600	0.142559388535373\\
700	0.120125347298123\\
800	0.127813088388748\\
900	0.101224401773305\\
1000	0.0935071652485435\\
};
\addlegendentry{QDEIM}
\addplot [color=mycolor2,line width=1.5pt, mark=asterisk, mark options={solid, mycolor2}]
  table[row sep=crcr]{%
500	0.13386588891567\\
600	0.108997653597618\\
700	0.0945115738426498\\
800	0.0837830141603245\\
900	0.0748905722434144\\
1000	0.0676767596289947\\
};
\addlegendentry{CADP-CUR}

\addplot [color=mycolor3,line width=1.5pt, mark=o, mark options={solid, mycolor3}]
  table[row sep=crcr]{%
500	0.135701321260413\\
600	0.108199161221196\\
700	0.092504604426317\\
800	0.0821069548885555\\
900	0.0757193503004416\\
1000	0.0687681378230886\\
};
\addlegendentry{DADP-CUR}

\addplot [color=mycolor4,line width=1.5pt, mark=square, mark options={solid, mycolor4}]
  table[row sep=crcr]{%
500	0.134833544959096\\
600	0.112387871545798\\
700	0.0953174801075694\\
800	0.0861417885274427\\
900	0.0782956399318748\\
1000	0.0675196615155772\\
};
\addlegendentry{CADP-CX}

\addplot [color=mycolor5,line width=1.5pt, mark=triangle, mark options={solid, rotate=180, mycolor5}]
  table[row sep=crcr]{%
500	0.132460301682578\\
600	0.107112820064372\\
700	0.0962099397252616\\
800	0.0836258815320712\\
900	0.0751315171815306\\
1000	0.0673570058562781\\
};
\addlegendentry{DADP-CX}

\end{groupplot}
\end{tikzpicture}%
}}
\caption{Relative approximation errors as a function of $k$ for the various iterative subselection DEIM CUR approximation algorithms compared with some standard CUR approximation algorithms using real data sets.\label{fig: 2}}
\end{figure}
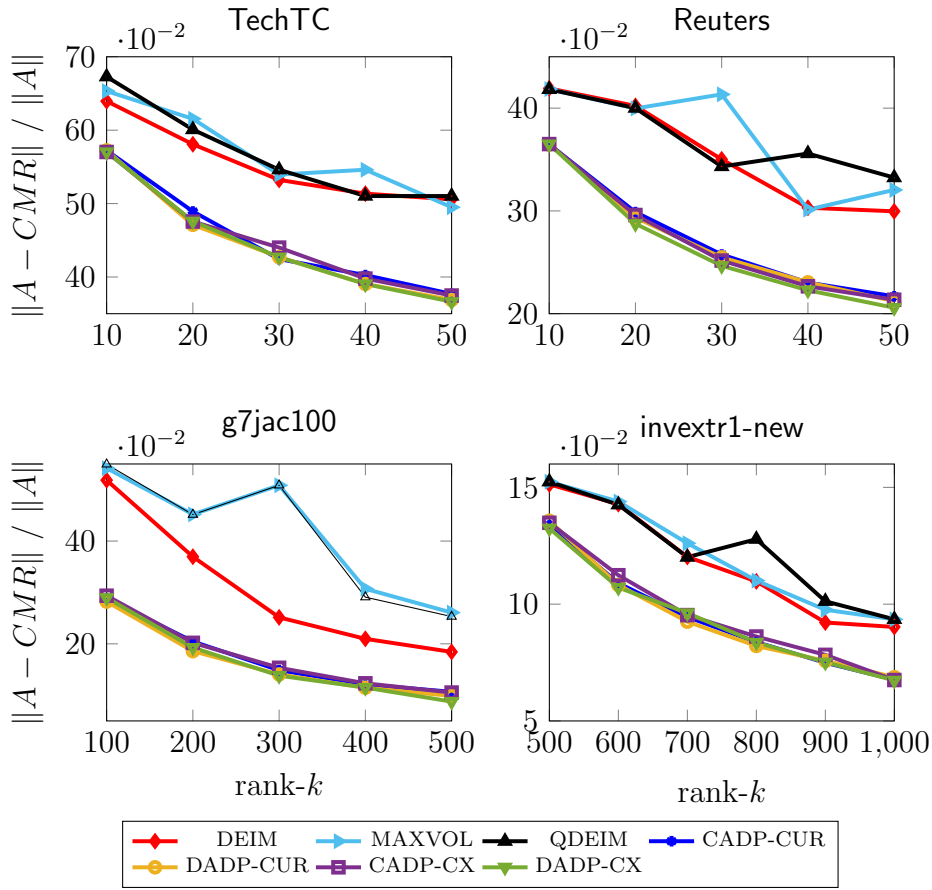

In \cref{fig: 3}, we maintain a constant rank while tweaking the user-defined parameters. To be more precise, for each data set we fix the rank and vary the values of the number of rounds ($t$) spanning from 2 to 6 for the \textbf{CADP-CX} and \textbf{CADP-CUR} algorithms, and $\delta$, which ranges from 0.4 to 0.8 for the \textbf{DADP-CX} and \textbf{DADP-CUR} methods. We observe that regardless of the specific values chosen for $t$ or $\delta$,  the key finding remains consistent. In all cases, the approximation errors of the proposed methods are lower than those associated with the standard one-round sampling methods.
\begin{figure}[htb!]
\centering
{{
%
%
\definecolor{mycolor1}{rgb}{0.00000,0.44700,0.74100}%
\definecolor{mycolor2}{rgb}{0,0,1}%
\definecolor{mycolor3}{rgb}{0.92900,0.69400,0.12500}%
\definecolor{mycolor4}{rgb}{0.49400,0.18400,0.55600}%
\definecolor{mycolor5}{rgb}{0.46600,0.67400,0.18800}%
\definecolor{mycolor6}{rgb}{0.30100,0.74500,0.93300}%
\definecolor{mycolor7}{rgb}{0,0,0}%
\begin{tikzpicture}

\begin{groupplot}[
group style={
    group name=my plots,
    group size=2 by 2,
    xlabels at=edge bottom,
    ylabels at=edge left,
    horizontal sep=1.3cm,vertical sep=3cm,
    },
    ymin=0.15,
    ymax=0.45,
    scaled y ticks=base 10:\exponent, 
    ylabel= {$\| A - CMR \|\ /\ \| A \|$},
legend style={at={(-0.2,-0.65)},anchor=south,legend  columns =4},
width=0.45\linewidth,
height=5cm
]
\nextgroupplot[title={\sf TechTC}, xlabel=rank-$30$,xmin=1,xmax=5,xtick={1,2,3,4,5},ymin=0.4,ymax=0.57]
\addplot [color=red, mark=diamond,line width=1.5pt,  mark options={solid, red}]
  table[row sep=crcr]{%
1	0.532279262004965\\
2	0.532279262004965\\
3	0.532279262004965\\
4	0.532279262004965\\
5	0.532279262004965\\
};

\addplot [color=mycolor6, mark=triangle,line width=1.5pt,  mark options={solid, rotate=270, mycolor6}]
  table[row sep=crcr]{%
1	0.539511434041749\\
2	0.539511434041749\\
3	0.539511434041749\\
4	0.539511434041749\\
5	0.539511434041749\\
};

\addplot [color=mycolor7, mark=triangle,line width=1.5pt,  mark options={solid, mycolor7}]
  table[row sep=crcr]{%
1	0.546136624140795\\
2	0.546136624140795\\
3	0.546136624140795\\
4	0.546136624140795\\
5	0.546136624140795\\
};

\addplot [color=mycolor2, mark=asterisk, line width=1.5pt, mark options={solid, mycolor2}]
  table[row sep=crcr]{%
1	0.475727114440007\\
2	0.439540653584266\\
3	0.439014895985175\\
4	0.443508913633526\\
5	0.440027481828186\\
};

\addplot [color=mycolor3, mark=o, line width=1.5pt, mark options={solid, mycolor3}]
  table[row sep=crcr]{%
1	0.451500005493798\\
2	0.447681817671309\\
3	0.439343123892232\\
4	0.441276980851511\\
5	0.426959321172227\\
};

\addplot [color=mycolor4, mark=square,line width=1.5pt,  mark options={solid, mycolor4}]
  table[row sep=crcr]{%
1	0.474220542942097\\
2	0.462956239230481\\
3	0.449392393997102\\
4	0.440736354211362\\
5	0.424326383229177\\
};

\addplot [color=mycolor5, mark=triangle,line width=1.5pt,  mark options={solid, rotate=180, mycolor5}]
  table[row sep=crcr]{%
1	0.470877088187415\\
2	0.456669604500147\\
3	0.450752953061929\\
4	0.437443985213455\\
5	0.427193142140243\\
};

\nextgroupplot[title={\sf Reuters}, xlabel=rank-$30$,xmin=1,xmax=5,xtick={1,2,3,4,5},ymin=0.2,ymax=0.45]
\addplot [color=red, mark=diamond, line width=1.5pt, mark options={solid, red}]
  table[row sep=crcr]{%
1	0.350182743753939\\
2	0.350182743753939\\
3	0.350182743753939\\
4	0.350182743753939\\
5	0.350182743753939\\
};

\addplot [color=mycolor6, mark=triangle,line width=1.5pt,  mark options={solid, rotate=270, mycolor6}]
  table[row sep=crcr]{%
1	0.41341658280015\\
2	0.41341658280015\\
3	0.41341658280015\\
4	0.41341658280015\\
5	0.41341658280015\\
};

\addplot [color=mycolor7, mark=triangle, line width=1.5pt, mark options={solid, mycolor7}]
  table[row sep=crcr]{%
1	0.343211415455115\\
2	0.343211415455115\\
3	0.343211415455115\\
4	0.343211415455115\\
5	0.343211415455115\\
};

\addplot [color=mycolor2, mark=asterisk,line width=1.5pt,  mark options={solid, mycolor2}]
  table[row sep=crcr]{%
1	0.296811656488167\\
2	0.271051593144557\\
3	0.272205745092776\\
4	0.267506151041207\\
5	0.251682756505136\\
};

\addplot [color=mycolor3, mark=o,line width=1.5pt,  mark options={solid, mycolor3}]
  table[row sep=crcr]{%
1	0.272135298129334\\
2	0.277799149287601\\
3	0.269421362338236\\
4	0.263566950899034\\
5	0.246599940620695\\
};

\addplot [color=mycolor4, mark=square,line width=1.5pt,  mark options={solid, mycolor4}]
  table[row sep=crcr]{%
1	0.325665499015754\\
2	0.277313014419358\\
3	0.281067429443239\\
4	0.281079044413708\\
5	0.257572601101225\\
};

\addplot [color=mycolor5, mark=triangle,line width=1.5pt,  mark options={solid, rotate=180, mycolor5}]
  table[row sep=crcr]{%
1	0.271242419075892\\
2	0.267206466213908\\
3	0.276157996226222\\
4	0.2685923286154\\
5	0.254779434003222\\
};

\nextgroupplot[title={\sf g7jac100}, xlabel=rank-$300$,ymin=0.05,ymax=0.6,xmin=1,xmax=5,xtick={1,2,3,4,5}]
\addplot [color=red, mark=diamond,line width=1.5pt,  mark options={solid, red}]
  table[row sep=crcr]{%
1	0.251480219806051\\
2	0.251480219806051\\
3	0.251480219806051\\
4	0.251480219806051\\
5	0.251480219806051\\
};

\addplot [color=mycolor6, mark=triangle,line width=1.5pt,  mark options={solid, rotate=270, mycolor6}]
  table[row sep=crcr]{%
1	0.508499732344036\\
2	0.508499732344036\\
3	0.508499732344036\\
4	0.508499732344036\\
5	0.508499732344036\\
};

\addplot [color=mycolor7, mark=triangle,line width=1.5pt,  mark options={solid, mycolor7}]
  table[row sep=crcr]{%
1	0.508501801217514\\
2	0.508501801217514\\
3	0.508501801217514\\
4	0.508501801217514\\
5	0.508501801217514\\
};

\addplot [color=mycolor2, mark=asterisk, line width=1.5pt, mark options={solid, mycolor2}]
  table[row sep=crcr]{%
1	0.239591323267642\\
2	0.217917006466928\\
3	0.179841995954225\\
4	0.180015852710137\\
5	0.15340342338346\\
};

\addplot [color=mycolor3, mark=o,line width=1.5pt,  mark options={solid, mycolor3}]
  table[row sep=crcr]{%
1	0.203732892741245\\
2	0.180031318871698\\
3	0.169565549759779\\
4	0.159591832701895\\
5	0.138052882601562\\
};

\addplot [color=mycolor4, mark=square,line width=1.5pt,  mark options={solid, mycolor4}]
  table[row sep=crcr]{%
1	0.228153526784649\\
2	0.218925565891686\\
3	0.186940088611982\\
4	0.17075445445595\\
5	0.148593211509307\\
};

\addplot [color=mycolor5, mark=triangle,line width=1.5pt,  mark options={solid, rotate=180, mycolor5}]
  table[row sep=crcr]{%
1	0.176947121048601\\
2	0.180812907178161\\
3	0.166213201529457\\
4	0.148585346253473\\
5	0.140477612418826\\
};

\nextgroupplot[title={\sf invextr1-new}, xlabel=rank-$600$,ymin=0.1,ymax=0.16,xmin=1,xmax=5,xtick={1,2,3,4,5}]
\addplot [color=red, mark=diamond, line width=1.5pt, mark options={solid, red}]
  table[row sep=crcr]{%
1	0.14282243851067\\
2	0.14282243851067\\
3	0.14282243851067\\
4	0.14282243851067\\
5	0.14282243851067\\
};
\addlegendentry{DEIM}

\addplot [color=mycolor6, mark=triangle,line width=1.5pt,  mark options={solid, rotate=270, mycolor6}]
  table[row sep=crcr]{%
1	0.143920494187863\\
2	0.143920494187863\\
3	0.143920494187863\\
4	0.143920494187863\\
5	0.143920494187863\\
};
\addlegendentry{MAXVOL}

\addplot [color=mycolor7, mark=triangle, line width=1.5pt, mark options={solid, mycolor7}]
  table[row sep=crcr]{%
1	0.142559388535373\\
2	0.142559388535373\\
3	0.142559388535373\\
4	0.142559388535373\\
5	0.142559388535373\\
};
\addlegendentry{QDEIM}

\addplot [color=mycolor2, mark=asterisk,line width=1.5pt,  mark options={solid, mycolor2}]
  table[row sep=crcr]{%
1	0.128349734780127\\
2	0.126098468835268\\
3	0.11342939411289\\
4	0.113377882661336\\
5	0.112387871545798\\
};
\addlegendentry{CADP-CX}

\addplot [color=mycolor3, mark=o, line width=1.5pt, mark options={solid, mycolor3}]
  table[row sep=crcr]{%
1	0.127114839057456\\
2	0.129275004064846\\
3	0.108947679980987\\
4	0.113234396421274\\
5	0.107112820064372\\
};
\addlegendentry{DADP-CX}

\addplot [color=mycolor4, mark=square,line width=1.5pt,  mark options={solid, mycolor4}]
  table[row sep=crcr]{%
1	0.130737001768937\\
2	0.12652245972584\\
3	0.112850153752989\\
4	0.111134376387531\\
5	0.108997653597618\\
};
\addlegendentry{CADP-CUR}

\addplot [color=mycolor5, mark=triangle,line width=1.5pt,  mark options={solid, rotate=180, mycolor5}]
  table[row sep=crcr]{%
1	0.131682230914091\\
2	0.117461875448869\\
3	0.109000037750735\\
4	0.113096969312227\\
5	0.108199161221196\\
};
\addlegendentry{DADP-CUR}

\end{groupplot}
\end{tikzpicture}%
}}

\caption{Relative approximation errors for the various iterative subselection DEIM CUR approximation algorithms compared with one-round sampling schemes
for a fixed rank $k$ with varying values of number of rounds $t=(2, 3, 5, 6, 10)$  for the \textbf{CADP-CX} and \textbf{CADP-CUR} algorithms and $\delta=(0.4,0.5,0.6,0.7,0.8)$  for the \textbf{DADP-CX} and \textbf{DADP-CUR} methods. 
\label{fig: 3}}
\end{figure}
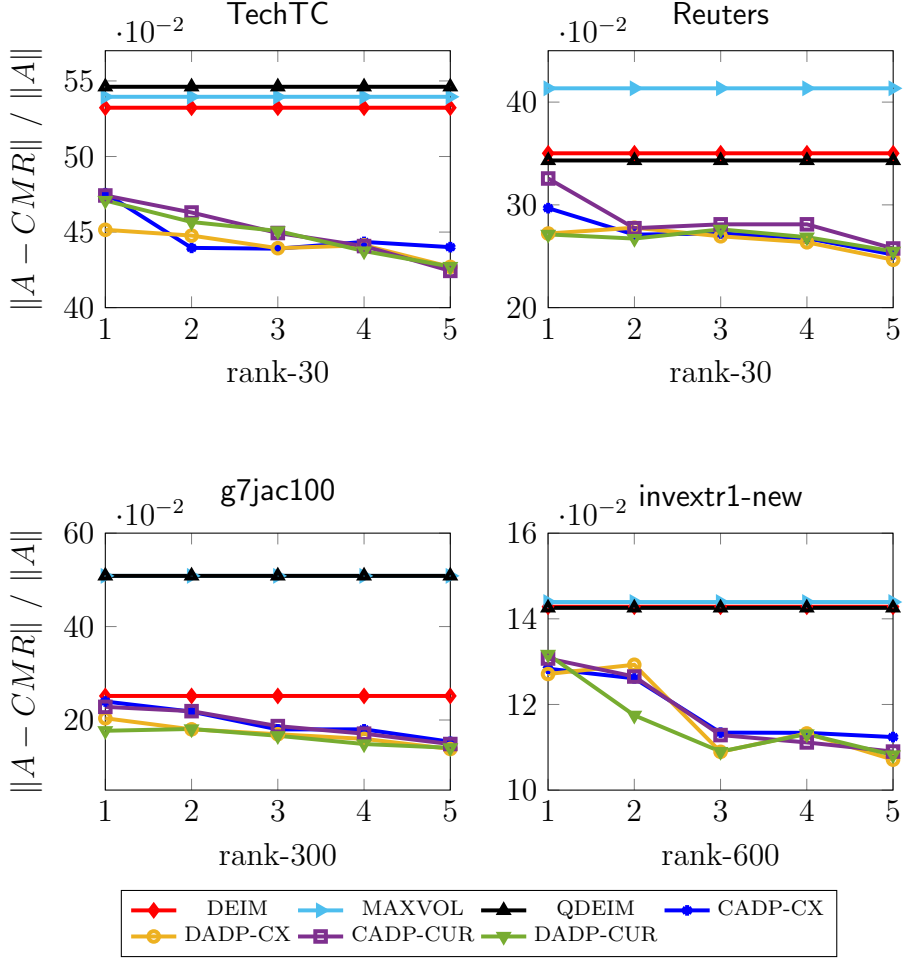

{\Cref{tab:exprand} presents a comparison between the randomized volume sampling \cref{algo:AdapSamp1} and deterministic iterative subselection algorithms proposed in this paper. Additionally, we include the results of a randomized variant of our proposed algorithm, termed CADP-CX-LVG, in which randomized leverage score sampling is employed instead of DEIM for index selection. To ensure fairness in comparison, we focus on algorithms requiring a constant number of indices per round, consistent with the volume sampling approach proposed in \cite{deshpande2006matrix}. The reported results pertain to the large-scale versions of CADP-CX and CADP-CUR (see, e.g., \cref{algo:NewAdapDEIM2}), wherein the residual is not explicitly computed. From the results, we observe that the deterministic algorithms yield lower approximation error while being computationally more efficient. Conversely, the volume sampling procedure produces higher approximation errors and is computationally more expensive since the residual $E=A-CC^+A$ is explicitly computed to update the probabilities.}

\begin{table}[htb!]
\centering
\caption{Relative approximation errors for the deterministic iterative subselection DEIM algorithms compared with randomized adaptive sampling methods. }
\label{tab:exprand}
{\footnotesize\begin{tabular}{lllllllll}\hline
Method/Data set   & \multicolumn{2}{c}{{\sf Reuters}} &  & \multicolumn{2}{c}{{\sf g7jac100}} &  & \multicolumn{2}{c}{{\sf invextr1-new}} \\
                  & \multicolumn{2}{c}{(rank-50)}     &  & \multicolumn{2}{c}{(rank-100)}     &  & \multicolumn{2}{c}{(rank-500)}         \\
                  &  Error            & Time (s)            &  &  Error            & Time(s)                &  &  Error              & Time (s)                 \\ \hline\rule{0pt}{2.5ex}%
Volume Sampling & 0.27                    & $3.64  \cdot 10^2$           &  & 0.66                  & $1.98\cdot 10^2$            &  & 0.32                    & $2.98 \cdot 10^3$              \\
CADP-CX-LVG       & 0.24                    & $1.70\cdot 10^1$              &  & 0.31                  & $7.54 \cdot 10^1$           &  & 0.18                    & $5.66 \cdot 10^2$              \\
CADP-CUR          & 0.17                   & $1.40\cdot 10^1$              &  & 0.29                  & $6.02 \cdot 10^1$            &  & 0.13                    & $4.50\cdot 10^2$              \\

CADP-CX           & 0.18                    & $1.80\cdot 10^1$              &  & 0.29                  & $7.01 \cdot 10^1$           &  & 0.13                    & $5.66 \cdot 10^2$ \\

\hline           
\end{tabular}}
\end{table}

As stated in \cref{sec:large-scale}, when dealing with large-scale matrices, performing a full (even reduced) SVD in each iteration of algorithms \ref{algo:AdapSamp}, \ref{algo:ExAdapDEIM}, and \ref{algo:NewAdapDEIM} can often become excessively costly. We evaluate the efficiency of the proposed algorithms: the small-scale versions, which compute the full SVD compared to their respective large-scale versions that use a Krylov--Schur routine to find a limited number of singular vectors (refer to \cref{algo:NewAdapDEIM2}). This experiment uses a large sparse data set, i.e.,~the {\sf Reuters-21578}. For our analysis, we approximate this matrix using a rank-50 approximation.

\cref{tab:smlascale} presents the results obtained from running the various algorithms. We observe that the large-scale variants (i.e., the various adaptations of \cref{algo:NewAdapDEIM2}), which utilize an iterative method for computing a few SVDs, demonstrate higher efficiency while maintaining similar approximation quality compared to the algorithms that compute the full SVD. Moreover, both for the full SVD and the iterative SVD routines, the algorithms with the residual defined as $E=A-CMR$ exhibit greater efficiency than those with the residual computed as $E=A-CC^+A$. Therefore, our new approach to computing the residual for the iterative subselection proves to be more efficient than the existing method while maintaining comparable approximation accuracy for this experiment.

\begin{table}[htb!]
\centering
\caption{Comparison of large-scale iterative subselection algorithms (iterative method for computing few SVDs) and small-scale iterative subselection algorithms (full SVD computation) on the {\sf Reuters-21578} data set approximation.}\label{tab:smlascale}
{\footnotesize
\begin{tabular}{l|ccc|cc}\hline \rule{0pt}{2.3ex}%
Method            & \multicolumn{2}{c}{Full SVD}       &  & \multicolumn{2}{c}{Iterative SVD}   \\ 
                  & Relative error & Runtime (s)       &  & Relative error & Runtime (s)        \\ \hline \rule{0pt}{2.3ex}%
\textbf{CADP-CUR} & $0.22$        & $1.18\cdot 10^3$ &  & $0.22$        & $1.34\cdot 10^1$  \\
\textbf{DADP-CUR} & $0.21$        & $2.36\cdot 10^3$ &  & $0.21$        & $1.31\cdot 10^1$  \\
\textbf{CADP-CX}  & $0.21$        & $2.27\cdot 10^3$ &  & $0.21$        & $1.63 \cdot 10^1$ \\
\textbf{DADP-CX}  & $0.21$        & $4.03\cdot 10^3$  &  & $0.21$        & $1.65\cdot 10^1$ \\\hline
\end{tabular}}
\end{table}

\noindent {\emph{\bf Interpretability of results.}
A primary objective of a CUR factorization is to construct factors that offer clear and interpretable insights. In this analysis, using the TechTC dataset, we compared the top 10 features selected by the DEIM algorithm with the iterative variants proposed in this paper. The aim is to evaluate whether terms selected by the algorithms are closely related to the categories of the documents, i.e., Evansville, Indiana (id 10567)
and Miami, Florida (id 11346). Here, the number of rounds is set to $k$ and $\delta=1$, implying that  CADP-CUR and DADP-CUR algorithms would yield identical results, as would CADP-CX and DADP-CX. \cref{tab:lvgscores} presents the leverage scores computed from the two leading singular vectors alongside the first 10 most significant columns (features) selected by the DEIM-based algorithms. The leading features indeed reveal key geographic terms closely related to the categories of the documents,  indicating that all algorithms select important features. While some terms are common across all approaches, the proposed algorithms tend to select terms with higher leverage scores compared to the one-round DEIM for terms that differ.

\begin{table}[htb!]
\centering
\caption{Comparison of the leading features selected by the DEIM-based algorithms for the TechTC dataset. The leverage scores (scaled) are computed using the leading two singular vectors.}\label{tab:lvgscores}
{\footnotesize\begin{tabular}{lllllllll}\hline \rule{0pt}{2.3ex}%
& \multicolumn{2}{c}{DEIM}    &  & \multicolumn{2}{c}{CADP-CX} &  & \multicolumn{2}{c}{CADP-CUR} \\
&Feature    & LVG Score &  & Feature          & LVG Score     &  & Feature          & LVG Score       \\  \hline \rule{0pt}{2.3ex}%
&evansville  & 1.000          &  & evansville       & 1.000              &  & evansville       & 1.000                \\
&florida     & 0.741          &  & florida          & 0.741              &  & florida          & 0.741                \\
&spacer      & 0.031          &  & spacer           & 0.031              &  & contact          & 0.055                \\
&contact     & 0.055          &  & about            & 0.073              &  & service          & 0.040                \\
&service     & 0.040          &  & information      & 0.113              &  & services         & 0.047                \\
&miami       & 0.058          &  & services         & 0.047              &  & please           & 0.015                \\
&chapter     & 0.004          &  & their            & 0.009              &  & spacer           & 0.031                \\
&health      & 0.005          &  & miami            & 0.058              &  & about            & 0.073                \\
&information & 0.113          &  & please           & 0.015              &  & indiana          & 0.030                \\
&events      & 0.013          &  & service          & 0.040              &  & south            & 0.123                  \\ \hline           
\end{tabular}}
\end{table}
}}
\end{experiment}

\section{Conclusions}
New approaches for selecting subsets of columns and rows using iterative subselection strategies have been presented. The first one is a DEIM adaptation of the so-called volume sampling \cite{deshpande2006matrix} for column subset selection. This procedure follows a fixed selection schedule, choosing a predetermined number of columns or rows in each iteration. In contrast, the second proposed iterative subselection strategy dynamically adjusts the selection schedule based on the decay of the singular values of the data. This approach aims to prioritize the selection of columns or rows that contribute the most to the overall structure and information of the data. By considering the significance of singular values and leveraging their decay pattern, the algorithm can adapt to the unique characteristics of the data, resulting in a more data-driven selection process.

Additionally, we also introduce an alternative approach for computing the residual in the index selection process. The first two iterative subselection algorithms we propose, i.e., \cref{algo:AdapSamp} and \cref{algo:ExAdapDEIM}, as well as existing adaptive sampling procedures \cite{optimalboutsidis,deshpande2006matrix,deshpande2006adaptive,paul2015column, Zhang}, define the residual as the error resulting from projecting the matrix $A$ onto either the column space of $C$ or the row space of $R$, i.e., $E = A - CC^+A$ or $E = A - AR^+R$, respectively. In contrast, our new method defines the residual as the error incurred by simultaneously projecting $A$ onto both the column space of $C$ and the row space of $R$. This entails computing a CUR factorization at each step using only the selected columns and rows.

We have also discussed how iterative procedures for computing a few singular vectors of large data matrices can be used with the newly proposed strategies. We have presented an adaptation of \cref{algo:ExAdapDEIM} for the large-scale case in \cref{algo:NewAdapDEIM2}, which can straightforwardly be adapted for \cref{algo:AdapSamp} and \cref{algo:NewAdapDEIM}. For each of the iterative subselection strategies proposed in this paper, we invoke the DEIM index selection method. However, we note that other deterministic index selection schemes such as the QDEIM technique \cite{Drmac} and the {\sf MaxVol} procedure \cite{Goreinov2010} may be employed. We have demonstrated through empirical analysis that the proposed methods in this work can produce better approximation results than the traditional method of \emph{one-round sampling} of all columns and rows. 

Overall, the proposed techniques may be useful for improving the accuracy of a CUR decomposition, but may also introduce additional complexities that need to be carefully addressed. The choice of whether to use the proposed iterative subselection methods or not may depend on the specific problem or application, as well as the trade-offs between accuracy, complexity, and computational resources.

 The data sets used in the experiments and the Matlab codes of the proposed algorithms are available via \href{https://https://github.com/perfectyayra/DEIM-with-Iterative-SVDs}{github.com/perfectyayra}.
\section*{Acknowledgements} 
This work has received funding from the European Union's Horizon 2020 research and innovation
programme under the Marie Sk\l odowska-Curie grant agreement No 812912.

\end{document}